\documentclass[reqno, 12pt, reqno]{amsart}

\usepackage{amsfonts, amsthm, amsmath, amssymb}
\usepackage{hyperref}
\hypersetup{colorlinks=false}

\usepackage[margin=1.5in]{geometry}

\RequirePackage{mathrsfs} \let\mathcal\mathscr

\numberwithin{equation}{section}
\newtheorem{theorem}{Theorem}[section]
\newtheorem{lemma}[theorem]{Lemma}

\theoremstyle{definition}

\newcommand{\card}{\#}
\newcommand{\ZZ}{\mathbb{Z}}

\newcommand{\NN}{\mathbb{N}}
\newcommand{\QQ}{\mathbb{Q}}

\newcommand{\ep}{\varepsilon}

\newcommand{\beql}[1]{\begin{equation}\label{#1}}
\newcommand{\eeq}{\end{equation}}

\begin{document}

\date{\today}

\title[Equidistribution for Solutions of $p+m^2+n^2=N$]
{Equidistribution for Solutions of $p+m^2+n^2=N$, and for
Ch\^{a}telet Surfaces} 
\author{D.R.\ Heath-Brown}

\address{Mathematical Institute\\
Radcliffe Observatory Quarter\\ Woodstock Road\\ Oxford\\ OX2 6GG}
\email{rhb@maths.ox.ac.uk}


\begin{flushright}
{\em For Henryk Iwaniec}\\
{\em In celebration of his 75th birthday}
\end{flushright}

\begin{abstract}

We present a general method for handling problems that ask for the
equidistribution of solutions to equations involving $m^2+n^2$, and
illustrate it by considering $p+m^2+n^2=N$.
\end{abstract}

\date{\today}
\maketitle

\section{Introduction}

It was shown by Linnik \cite{Linnik}, \cite{Linnik2} that every large
integer $N$ may
be represented as $N=p+m^2+n^2$, with $p$ prime. His work built on a
method of Hooley \cite{Hooley}, who had established this result
subject to the Generalized Riemann Hypothesis. In fact there is a
positive constant $\delta$ such that 
\[\card\{(p,m,n):p+m^2+n^2=N\}=\frac{N}{\log N}C(N)
+O(N(\log N)^{-1-\delta})\]
with
\[C(N)=\pi\prod_{p>2}\left(1+\frac{\chi_4(p)}{p(p-1)}\right)
\prod_{p\mid N}\frac{(p-1)(p-\chi_4(p))}{p^2-p+\chi_4(p)},\]
(where $\chi_4$ is the non-trivial character modulo 4) --- 
see Hooley \cite[Chapter 5]{Hooley2}, for example. We remark here
that
\[(\log\log N)^{-1/2}\ll C(N)\ll (\log\log N)^{1/2},\]
so that the asymptotic formula saves a power of $\log N$ over the main
term. 

In this paper we
examine the distribution of the solutions $p,m,n$. There are various
ways to do this, but the approach we describe appears to be very
flexible. It should be stressed however that most of the results we
describe in the context of the equation $N=p+m^2+n^2$ and its
generalizations can be obtained by a different route, given by
Bredikhin and Linnik \cite{BL}.

There is no difficulty in controlling the size of the variable $p$, and
it is an easy matter to adapt Hooley's argument to show that
\[\card\{(p,m,n):p+m^2+n^2=N,\, aN<p\le bN\}\hspace{3cm}\]
\[\hspace{3cm} {}=(b-a)\frac{N}{\log N}C(N)
+O(N(\log N)^{-1-\delta})\]
uniformly for $0\le a<b\le 1$. If $aN<p\le bN$ then
$(m,n)$ will lie in the annulus $(1-b)N\le ||(m,n)||_2\le (1-a)N$. In
order to describe the equidistribution of $(m,n)$ within this annulus
it is natural to write $\alpha=m+in\in\ZZ[i]$, with
$N(\alpha)=m^2+n^2$, and to consider $\arg(\alpha)$. We then have the
following result.
\begin{theorem}\label{th1}
  There is an absolute constant $\delta>0$ as follows.
  Let real numbers $a,b,c,d$ be given, with $0\le a<b\le 1$ and
  $0\le c<d\le 2\pi$. Let $N\in\NN$ be given. Then 
  \[\card\{(p,\alpha):\, p \mbox{ prime},\,\alpha\in\ZZ[i],\,
  p+N(\alpha)=N,\, aN<p\le bN,\, c<\arg(\alpha)\le d\}\]
  \[\hspace{2cm} {}=
  \frac{(b-a)(d-c)}{2\pi}\frac{N}{\log N}C(N)
+O_{\ep}(N(\log N)^{-1-\delta}).\]
\end{theorem}

The natural way to handle the restriction on $\arg(\alpha)$ is via the
Erd\H{o}s--Tur\'{a}n inequality. For any finite subset $A$ of $\ZZ[i]$
and any $H\in\NN$ we have
\[\card\{\alpha\in A:\, c<\arg(\alpha)\le d\}\hspace{6cm}\]
\[\hspace{3cm} {}=\frac{d-c}{2\pi}\card A
+O(H^{-1}\card A)
+O\left(H^{-1}\sum_{h\le H}\left|\sum_{\alpha\in A}
\left(\frac{\alpha}{|\alpha|}\right)^h\right|\right).\]
In our setting the term $O(H^{-1}\card A)$ will contribute
$O(NC(N)/(H\log N))$ in Theorem \ref{th1}, and so we see that the required
result will follow from the following estimate, on taking $H=\log N$, for example.
\begin{theorem}\label{th2}
Define
\[g_h(n)=\frac{1}{4}\sum_{N(\alpha)=n}
\left(\frac{\alpha}{|\alpha|}\right)^h\]
and $f_h(n)=|g_h(n)|$. Then
  \[\sum_{p<N}f_h(N-p)  \ll_{\ep} N(\log N)^{-5/4+\ep}\]
for any fixed $\ep>0$,  uniformly for $1\le h\le\log N$.
\end{theorem}
We remark that the exponent $5/4$ can be improved to $2-2/\pi$, but
this has no qualitative benefit for us. 

The reader may be surprised that it is sufficient to consider the average of 
$f_h(N-p)$, rather than $g_h(N-p)$. The fact that
this is indeed enough is the fundamental observation behind 
this paper. As motivation we point out that for
most primes $p$ the number $N-p$ will not be a sum of two squares, and
$g_h(N-p)$ will therefore vanish. However when $N-p$ is a sum of
two squares the sum for $g_h(N-p)$ will typically contain many terms
$\alpha$, so that we may expect some cancellation to occur.

The above idea, that we can obtain enough cancellation by considering
$|g_h(n)|$ rather than $g_h(n)$, has occurred previously in other contexts.
(We are grateful to Valentin Blomer for this observation.)
Thus, for example, Holowinsky's work \cite{holo} on Quantum Unique 
Ergodicity looks at shifted convolution sums
\[\sum_{n\le x}|\lambda_1(n)\lambda_2(n+l)|,\]
which can be handled non-trivially using Erd\H{o}s' method. In 
Holowinsky's application one takes $\lambda_1=\lambda_2$ to be
the Fourier coefficients of a holomorphic modular cusp form.

\bigskip

The above ideas can be applied in many other situations. For example,
one may consider the equidistribution of points on Ch\^{a}telet
surfaces, given by equations
\[F(u,v)=x^2-ay^2,\]
where $F\in\ZZ[u,v]$ is a separable quartic polynomial and $a\in\ZZ$ is
not a square. Investigations into Manin's Conjecture for these
surfaces have focused on the case $a=-1$, but in this situation the
conjecture has been established for all forms $F$. The most difficult
case, in which $F$ is irreducible over $\QQ(i)$, has been handled by 
de la Bret\`{e}che and Tenenbaum \cite{dlBT}. Since one is interested
in rational points it is natural to consider solutions to
\beql{c1}
t^2F(u,v)=x^2+y^2
\eeq
with
\beql{c2}
t,u,v,x,y \in\ZZ,\;\;\; \gcd(u,v)=\gcd(t,x,y)=1,\;\;\; t>0.
\eeq
The counting function considered by de la Bret\`{e}che and Tenenbaum,
is one half the number of solutions satisfying the height condition
\beql{hb}
t^{1/2}\max(|u|,|v|)\le B^{1/2}.
\eeq
They then obtain an asymptotic formula
\[\sigma_{\infty}(F)\mathfrak{S}(F)B\log B+O_F(B(\log B)^{99/100}),\]
where
\[\sigma_{\infty}(F)=
\frac{\pi}{2}\mathrm{meas}\{(u,v)\in [-1,1]^2:\, F(u,v)>0\}\]
and $\mathfrak{S}(F)$ is a product of $p$-adic densities. There is no
difficulty in adapting the argument to replace (\ref{hb}) with a
condition
\beql{reg}
t^{1/2}(u,v)\in B^{1/2}\mathcal{R}
\eeq
for an arbitrary rectangle $\mathcal{R}\subseteq[-1,1]^2$. This
produces an analogous asymptotic formula with $\mathfrak{S}(F)$
replaced by
\[\sigma_{\infty}(F;\mathcal{R})=
\frac{\pi}{2}\mathrm{meas}\{(u,v)\in\mathcal{R}:\, F(u,v)>0\}.\]
We then have the following result, which controls the argument of
$x+iy$, in addition to the location of $(u,v)$.
\begin{theorem}\label{th3}
  Suppose one is given a rectangle $\mathcal{R}\subseteq[-1,1]^2$ and
  real numbers $a<b$ in $[0,2\pi]$. Let $N(B)$ be half the number of
    5-tuples $(t,u,v,x,y)\in\ZZ^5$ satisfying (\ref{c1}) and
    (\ref{c2}), and with
    \[t^{1/2}(u,v)\in  B^{1/2}\mathcal{R}\;\;\;\mbox{and}\;\;\;
    \arg(x+iy)\in[a,b].\]
    Then
    \[N(B)=\frac{b-a}{2\pi}
\sigma_{\infty}(F;\mathcal{R})\mathfrak{S}(F)B\log B+O_F(B(\log B)^{99/100}).\]
\end{theorem}
The key input here is the following estimate, which we will prove in
\S \ref{Pth2}
\begin{lemma}\label{lt2}
  For any fixed $\ep>0$ we have
  \[\sum_{t^{1/2}\max(|u|,|v|)\le B^{1/2}}
  \left|\sum_{\substack{x^2+y^2=t^2F(u,v)\\ \gcd(t,x,y)=1}}
  \left(\frac{x+iy}{|x+iy|}\right)^{4h}\right|\ll_{F,\ep}
  B(\log B)^{3/4+\ep}\]
    uniformly for $1\le h\le\log B$.
\end{lemma}

These ideas are not restricted to the use of
Gr\"{o}ssencharacters. For example we can consider the equation
$N=p+m^2+n^2$ subject to congruence conditions on $m$ and $n$ by
looking at $N=p+N(\alpha)$ subject to $\alpha\equiv\alpha_0($mod $k)$,
say, with $\alpha_0\in\ZZ[i]$ coprime to $k\in\NN$. This condition can
be detected with Dirichlet characters $\chi($mod $k)$ over
$\ZZ[i]$. With these we have the following estimate.
\begin{theorem}\label{t2c}
  Let $k\in\NN$ and $\ep>0$ be given.  Then
\[\sum_{p<N}\left|\frac{1}{4}\sum_{N(\alpha)=N-p}\chi(\alpha)\right|
\ll_{k,\ep} N(\log N)^{-5/4+\ep}\]
for every non-principal Dirichlet character $\chi($mod $k)$
over $\ZZ[i]$ such that $\chi(i)=1$.
\end{theorem}
We leave the details of the proof to the reader. In order to apply
this to examine $N=p+N(\alpha)$ with $\alpha\equiv\alpha_0($mod $k)$,
(for $\alpha_0$ coprime to $k$) we need to know about
\[\sum_{N=p+N(\alpha)}\chi_0(\alpha)\]
where $\chi_0$ is the principal character modulo $k$. However $\alpha$
and $k$ are coprime if and only if $N(\alpha)$ and $k$ are coprime, so
that it is enough to consider
\[\card\{(p,m,n):\, N=p+m^2+n^2,\, \gcd(N-p,k)=1\}.\]
This may be tackled by a trivial variant of Hooley's method, in which
one puts congruence conditions on $p$.

One can also handle congruence conditions $\alpha\equiv\alpha_0($mod
$k)$ when $\alpha$ shares a common factor with $k$. If
$\alpha_0=\beta\alpha_1$ and $k=\beta\gamma$ with $\alpha_1$ and
$\gamma$ coprime then the congruence $\alpha\equiv\alpha_0($mod
$k)$ implies that $\alpha=\beta\alpha'$, say, with
\beql{cc}
\alpha'\equiv\alpha_1(\mathrm{mod}\,\gamma).
\eeq
The equation $N=p+N(\alpha)$
then requires $p\equiv N($mod $N(\beta))$, so that we have to solve
\[\frac{N-p}{N(\beta)}=N(\alpha')\;\;\;\mbox{with}
\;\;\; p\equiv N(\mathrm{mod}\, N(\beta)),\]
subject to the congruence condition (\ref{cc}). We can tackle this
with the machinery described above, using a version of Theorem
\ref{t2c} in which one restricts to primes $p\equiv N($mod $N(\beta))$
and replaces the summation condition $N(\alpha)=N-p$ by
$N(\alpha)=(N-p)/N(\beta)$. In particular this is enough to give an
asymptotic formula for the number of solutions to
\[N=p+a^4m^2+b^4n^2\]
for any fixed $a,b\in\NN$, and thereby to count solutions of
\[N=p+m^2+n^2,\;\;\;\; m,n\mbox{ square-free}.\]
This is the problem for which Hooley \cite{hooley3} proves a positive
lower bound. The reader should note though that Hooley
\cite[page 202]{hooley3} suggests that the necessary asymptotic
formulae ``can probably be established \ldots by a different method
involving an elaborate use of the arithmetic of binary quadratic forms.''

One can also replace $m^2+n^2$ by other quadratic forms, and this was
the main achievement (in the context of the equation $N=p+m^2+n^2$)
in the work of Bredikhin and Linnik \cite{BL}
mentioned earlier. As an example, we consider the Ch\^{a}telet surface
\[F(u,v)=x^2-2y^2,\]
where $F(u,v)\in\ZZ[u,v]$ is a quartic form, assumed now to be irreducible
over $\QQ(\sqrt{2})$. The function
\[r_2(n)=\sum_{\substack{d\mid n\\ d\,\mathrm{odd}}}\left(\frac{2}{d}\right),
\;\;\; n\in\ZZ\setminus\{0\},\]
counts solutions of $n=x^2-2y^2$ once from each set of associates of
$x+y\sqrt{2}$. Equivalently, $r_2(n)$ counts integral ideals $I$ of
$\QQ(\sqrt{2})$ having norm $N(I)=|n|$. (Previously we had used $N(*)$
for the norm function on $\QQ(i)$, but here it is the norm on
$\QQ(\sqrt{2})$. We trust this will cause no confusion.)
A moment's reflection shows that if $t^2\mid x^2-2y^2$ then
$\gcd(t,x,y)=1$ if and only if $\gcd(t,\alpha,\alpha^{\sigma})=1$,
where $\alpha=x+y\sqrt{2}$ and $\sigma$ is the non-trivial
automorphism of $\QQ(\sqrt{2})$.
One therefore wants to count solutions to
\beql{Ie}
t^2|F(u,v)|=N(I)
\eeq
subject to $\gcd(u,v)=\gcd(t,I,I^{\sigma})=1$, lying in a region (\ref{reg}),
these being the natural analogues of (\ref{c1}), (\ref{c2}).
The method of de la Bret\`{e}che and Tenenbaum
\cite{dlBT} can be readily adapted to show that the number of
solutions takes the form
\[c(\mathcal{R})B\log B+O_F(B(\log B)^{99/100}),\]
with a constant $c(\mathcal{R})$ depending on $F$ as well as $\mathcal{R}$.
However one would wish to control the location of the corresponding
points $(x,y)$ for which $I=(x+y\sqrt{2})$. If $\mathcal{R}$ is
suitably small one will know the product
\[\left(\frac{|x-y\sqrt{2}|}{t}\right)\left(\frac{|x+y\sqrt{2}|}{t}\right)
=t^{-2}N(I)=|F(u,v)|\]
to a good degree of approximation, and hence it is enough to control
the quotient
\[\left.\left(\frac{|x-y\sqrt{2}|}{t}\right)\middle/
\left(\frac{|x+y\sqrt{2}|}{t}\right)=
\frac{|x-y\sqrt{2}|}{|x+y\sqrt{2}|}.\right.\]
We can do this, and indeed slightly more, by using the
Gr\"{o}ssencharacter
\[\chi(I)=\mathrm{sgn}(x^2-2y^2)
\exp\left\{\pi
i\frac{\log|x-y\sqrt{2}|-\log|x+y\sqrt{2}|}{2\log(1+\sqrt{2})}\right\},\]
where $I=(x+y\sqrt{2})$.  (The reader who is unfamiliar with such things 
may readily check that this is at least well-defined.)
In analogy to Lemma \ref{lt2} we have the following bound.
\begin{lemma}\label{lt?}
  For any fixed $\ep>0$ we have
  \[\sum_{t^{1/2}\max(|u|,|v|)\le B^{1/2}}
  \left|\sum_{\substack{N(I)=t^2F(u,v)\\ \gcd(t,I,I^{\sigma})=1}}
  \chi(I)^{h}\right|\ll_{F,\ep}
  B(\log B)^{3/4+\ep}\]
    uniformly for $1\le h\le\log B$.
\end{lemma}
This allows us to count solutions to (\ref{Ie}) asymptotically, with
$\chi(I)$ restricted to any given arc of the unit circle. A moment's
thought shows that this corresponds to having $I=(x+y\sqrt{2})$ with
$(x-y\sqrt{2})/(x+y\sqrt{2})$ in a given interval
$(a,b)\subset(-\infty,0)\cup(0,\infty)$. We are then able to count
solutions to $t^2F(u,v)=x^2-2y^2$ in which both $(u,v)$ and $(x,y)$
lie in prescribed regions.

Lastly we mention that one can work analogously with class group
characters. Suppose, for example, that we are interested in the
Ch\^{a}telet surface
\[F(u,v)=x^2+14y^2,\]
for a quartic form $F\in\ZZ[u,v]$ which is irreducible over
$\QQ(\sqrt{-14})$. \footnote{We are grateful to Evan O'Dorney for
  pointing out an error in the discussion of this problem in an
  earlier version of this paper.}
There are four classes of positive definite binary
quadratic forms of discriminant $-56$, with representatives $x^2+14y^2$,
$2x^2+7y^2$, $3x^2+2xy+5y^2$,
and $3x^2-2xy+5y^2$. Since the first two of these lie in the same genus
it is not
possible to separate out representations by these two forms using only
congruence conditions. The combined number of representations of $n\in\NN$
by all four classes of forms is twice the number of
integral ideals $I$ of $\QQ(\sqrt{-14})$ with norm $N(I)=n$; and the
number of representations by $x^2+14y^2$ is twice the number of
such $I$ that are principal. As above we will want to consider
solutions to $t^2F(u,v)=N(I)$ with $\gcd(u,v)=1$ and
$\gcd(t,I,I^{\sigma})=1$, where $\sigma$ is now complex
conjugation. We may obtain an asymptotic formula for the problem in
which we count both principal and non-principal ideals $I$, based on
the formula
\[\card\{I:\, N(I)=n\}=\sum_{\substack{d\mid n\\ \gcd(d,14)=1}}
\left(\frac{-14}{d}\right),\]
and following the method laid out by de la Bret\`{e}che and Tenenbaum
\cite{dlBT}. Then, to restrict the count to principal ideals we use
class group characters.  Following
our previous argument we then see that the following lemma suffices
\begin{lemma}\label{cgc}
Let $\chi$ be a non-trivial class group character for
$\QQ(\sqrt{-14})$, and let $N$ be the norm function for
$\QQ(\sqrt{-14})$. Then for any fixed $\ep>0$ we have
  \[\sum_{t^{1/2}\max(|u|,|v|)\le B^{1/2}}
  \left|\sum_{\substack{N(I)=t^2F(u,v)\\ \gcd(t,I,I^{\sigma})=1}}
  \chi(I)\right|\ll_{F,\ep}
  B(\log B)^{3/4+\ep},\]
where $I$ runs over integral ideals of $\QQ(\sqrt{-14})$.
\end{lemma}
This may be proved following the same reasoning as for Lemma
\ref{lt?}, the key information input being the Prime Ideal Theorem for
ideal classes.

\section{Proof of Theorem \ref{th2}}\label{pt2}

It is clear that $g_h(n)$ vanishes unless $4\mid h$, since each
$\alpha$ has associates $\alpha, i\alpha, -\alpha, -i\alpha$, and
$1+i^h+(-1)^h+(-i)^h=0$ unless $4\mid h$. We therefore assume that
$4\mid h$ henceforth. We now claim that $g_h(n)$
is multiplicative, whence $f_h(n)$ is also multiplicative. To
establish the claim, suppose that $n=uv$ with coprime factors
$u,v$. Then it merely suffices to observe that
each $\alpha$ of norm $n$ can be written as
$\alpha=\beta\gamma$ with $N(\beta)=u$ and $N(\gamma)=v$ in exactly four
ways, corresponding to the four associate choices of
$\beta=\mathrm{g.c.d.}(\alpha,u)$. (We leave it to the reader to check
this.)

We now call on a general result that gives accurate order of magnitude
estimates for sums of non-negative multiplicative functions. Such
bounds go back to Erd\H{o}s \cite{erdos}, and were investigated
further by Barban and Vehov \cite{BV}, and Shiu \cite{shiu} amongst
others. For our application we need a version that applies to sums of
functions over sifted sequences and this has been given by 
Pollack \cite[Theorem 1.1]{pollack}.  In his notation we take $k=1$,
$\beta=\tfrac14$ and $y=x$.
\begin{lemma}
  Let $f(n)$ be a multiplicative function satisfying
  $0\le f(n)\le d(n)$ for all $n\in\NN$. Suppose that, for each prime
  $p\le x$ the set $\mathcal{E}_p$ is either the empty set or is a
  non-zero residue class modulo $p$, and write $\nu(p)=0$ or
  $\nu(p)=1$ accordingly.  Let
  \[\mathcal{S}=\bigcap_{p\le x}\mathcal{E}_p^c.\]
  Then
  \[\sum_{\substack{n\le x\\ n\in\mathcal{S}}}f(n)\ll
  \frac{x}{\log x}\exp\left(\sum_{p\le x}\frac{f(p)-\nu(p)}{p}\right).\]
\end{lemma}
The reader should observe that our requirement that $0\le f(n)\le d(n)$
implies Pollack's condition $f\in\mathcal{M}$.

We apply this lemma with $x=N$ and $f=f_h$.  We take $\mathcal{E}_p$
to consist
of the residue class $N$(mod $p$) when $p\le\sqrt{N}$ does not
divide $N$, and to be the empty set otherwise. Then if $p\ge\sqrt{N}$
is prime we will have $N-p\in\mathcal{S}$. The lemma then tells us
that
\[\sum_{\sqrt{N}\le p\le N}f_h(N-p)\ll
  \frac{N}{\log N}\exp\left(\sum_{p\le N}\frac{f_h(p)-\nu(p)}{p}\right).\]
    Here we note that
  \begin{eqnarray}\label{lll}
    \sum_{p\le N}p^{-1}\nu(p)&\ge&\sum_{p\le\sqrt{N}}p^{-1}-
    \sum_{p\mid N}p^{-1}\nonumber\\
    &\ge&\log\log N+O(1)-\sum_{p\le\omega(N)}p^{-1}\nonumber\\
    &\ge&\log\log N+O(1)-\sum_{p\le\log N}p^{-1}\nonumber\\
    &=&\log\log N-\log\log\log N+O(1).
  \end{eqnarray}
  so that
\[\sum_{\sqrt{N}\le p\le N}f_h(N-p)\ll \frac{N\log\log N}{(\log N)^2}
\exp\left(\sum_{p\le N}\frac{f_h(p)}{p}\right)\]
and hence
\beql{e2}\sum_{p\le N}f_h(N-p)\ll\sqrt{N}+ \frac{N\log\log N}{(\log N)^2}
\exp\left(\sum_{p\le N}\frac{f_h(p)}{p}\right).
\eeq
It remains to consider
\[\sum_{p\le N}\frac{f_h(p)}{p}.\]
Of course $f_h(p)=0$ if $p\equiv 3($mod 4$)$, while if
$p\equiv 1($mod 4$)$ we may write $p=N(\alpha)$ for some particular
$\alpha$, in which case $f_h(p)=2|\cos(h\arg(\alpha))|$.  We now use
the famous inequality $3+4\cos\theta+\cos 2\theta\ge 0$, which implies
\[2|\cos\theta|\le\tfrac32+\tfrac12\cos(2\theta)\]
for all real $\theta$, and hence
\beql{34}
f_h(p)\le\tfrac32+\tfrac14 g_{2h}(p),
\eeq
for $p\equiv 1($mod 4$)$. Since $f_h(p)$ and $g_{2h}(p)$ vanish for 
$p\equiv 3($mod 4$)$ it follows that
\[\sum_{3\le p\le N}\frac{f_h(p)}{p}\le
\frac32\sum_{\substack{p\le N\\p\equiv 1(\mathrm{mod}\, 4)}}\frac{1}{p}
+\frac14\sum_{3\le p\le N}\frac{g_{2h}(p)}{p}.\]
The first sum on
the right is $\tfrac12 \log\log N+O(1)$, whence
\beql{e1}
\sum_{3\le p\le N}\frac{f_h(p)}{p}\le
\frac34\log\log N+O(1)+\frac14\sum_{3\le p\le N}\frac{g_{2h}(p)}{p}.
\eeq
To estimate the sum on the right we
call on the Prime Number Theorem for primes over $\QQ(i)$
with Gr\"{o}ssencharacter, in the following form.
\begin{lemma}\label{pntg}
  There is an absolute constant $c>0$ such that
  \[\sum_{\substack{\alpha\in\ZZ[i]\\ N(\alpha)\le x}}\Lambda(\alpha)
  \left(\frac{\alpha}{|\alpha|}\right)^{4k}\ll
  x\exp\left\{-c\frac{\log x}{\sqrt{\log x}+\log k}\right\}(\log xk)^4\]
  uniformly for $x\ge 2$ and $k\in\NN$.
  \end{lemma}
 This is a special case of Theorem 5.13 in Iwaniec and Kowalski
 \cite{IK}, as the reader may confirm.  Using partial summation we
 may deduce from Lemma
 \ref{pntg} that if $N_0=\exp\{(\log 3k)^2\}$ then
\[\sum_{\substack{\alpha\in\ZZ[i]\\ N_0<N(\alpha)\le N}}
\frac{\Lambda(\alpha)}{N(\alpha)\log N(\alpha)}
\left(\frac{\alpha}{|\alpha|}\right)^{4k}\ll 1\]
uniformly for $N\ge N_0$ and $k\in\NN$. A trivial bound shows that
we may include terms with $N(\alpha)\le N_0$ at a cost
$O(\log\log k)$, and that we may remove terms
in which $N(\alpha)$ is not a prime at a cost $O(1)$. Taking $4k=2h$
we then deduce that
\[\sum_{3\le p\le N}\frac{g_{2h}(p)}{p}\ll \log\log h\]
whenever $4\mid h$. Substituting this into (\ref{e1}) we find that
\[\sum_{3\le p\le N}\frac{f_h(p)}{p}\le
\frac{3}{4}\log\log N+O(\log\log h).\]
Inserting this bound into (\ref{e2}) we have
\[\sum_{p\le N}f_h(N-p)\ll\sqrt{N}+ \frac{N\log\log N}{(\log N)^{5/4}}(\log h)^A\]
for some absolute constant $A$. This is sufficient for Theorem \ref{th2}.

\section{Proof of Theorem \ref{th3}}\label{Pth2}

In analogy to the argument that leads from Theorem \ref{th2} to
Theorem \ref{th1} we see that Lemma \ref{lt2} will suffice for
the proof of Theorem \ref{th3}. To establish the lemma we first
note that $t$ is
composed entirely of primes $p\equiv 1($mod 4) whenever
$t^2\mid x^2+y^2$ with $\gcd(t,x,y)=1$.  Moreover if
$p^e\mid t$ with a prime $p=N(\pi)$ over $\ZZ[i]$, then we must have
either $\pi^e\mid x+iy$ or $\pi^e\mid x-iy$. Furthermore, if
$\pi^e\mid x+iy$ with $e\ge 1$ it follows that $x+iy$ must be coprime
to $\overline{\pi}$. Thus $x+iy$ must be divisible by some $\mu$ with
$N(\mu)=t$, such that $\gcd(\mu,\overline{\mu})=1$, and with
$\gcd(x+iy,\overline{\mu})=1$. In order to count each set of
associates $i^n\mu$ just once we restrict $\mu$ to be primary; that is
to say we require that $\mu\equiv 1($mod $2+2i)$.

We now write $x+iy=\mu\alpha$, so that
\[\sum_{\substack{x^2+y^2=t^2F(u,v)\\ \gcd(t,x,y)=1}}
  \left(\frac{x+iy}{|x+iy|}\right)^{4h}=
  \sum_{\mu}\left(\frac{\mu}{|\mu|}\right)^{4h}
\sum_{\substack{N(\alpha)=F(u,v)\\ \gcd(\alpha,\overline{\mu})=1}}
\left(\frac{\alpha}{|\alpha|}\right)^{4h},\]
the sum over $\mu$ being for solutions of $N(\mu)=t$ subject to
$\gcd(\mu,\overline{\mu})=1$ and $\mu\equiv 1($mod $2+2i)$. We
therefore set
\[f_k(n;\mu)=\frac14
\left|\sum_{\substack{N(\alpha)=n\\ \gcd(\alpha,\overline{\mu})=1}}
\left(\frac{\alpha}{|\alpha|}\right)^k\right|,\]
which differs from our previous function $f_h(n)$ only in the
condition $\gcd(\alpha,\overline{\mu})=1$. It now follows that
\[\sum_{t^{1/2}\max(|u|,|v|)\le B^{1/2}}
  \left|\sum_{\substack{x^2+y^2=t^2F(u,v)\\ \gcd(t,x,y)=1}}
    \left(\frac{x+iy}{|x+iy|}\right)^{4h}\right|\hspace{2cm}\]
    \beql{yi}
    \hspace{2cm}{}\le 4\sum_{\substack{\gcd(\mu,\overline{\mu})=1\\ \mu\equiv
    1(\mathrm{mod}\, 2+2i)}}\;\;\sum_{\max(|u|,|v|)\le
      (B/N(\mu))^{1/2}}f_{4h}(F(u,v);\mu).
    \eeq
The inner sum above can now be bounded by an Erd\H{o}s type estimate,
a suitable version of which is given in the next lemma, which is an
immediate consequence of de la Bret\`{e}che and Browning
\cite[Corollary 1]{BdlB}.
\begin{lemma}\label{dlBB}
  Let $F(u,v)\in\ZZ[u,v]$ be an irreducible binary form of degree 
  $d\ge 2$, and let $f(n)$ be a multiplicative function satisfying
  $0\le f(n)\le d(n)$ for all $n\in\NN$.  Then
  \[\sum_{\max(|u|,|v|)\le X}f(|F(u,v)|)\ll_F
  X^2\prod_{d<p\le X}\left(1+\frac{\rho(p)(f(p)-1)}{p}\right),\]
  where
  \[\rho(p)=\frac{1}{p-1}\card\{(u,v)\in(0,p]^2:\, p\mid F(u,v),\,
    \gcd(u,v,p)=1\}.\]
\end{lemma}

The reader should note that we have $d_1=d_2=0$ and $F=G$ in the notation of
\cite{BdlB}. We should also stress that the implied constant in Lemma
\ref{dlBB} depends only on $F$, and is uniform over all function $f$
satisfying $0\le f(n)\le d(n)$.

In our situation $\rho(p)$ is the number of zeros modulo $p$ of the
polynomial $F(X,1)$ as soon as $p>|F(1,0)|$. Since $\rho(p)$ and
$f_{4h}(*;\mu)$ are absolutely bounded we see that Lemma \ref{dlBB}
yields
\[\sum_{\max(|u|,|v|)\le (B/N(\mu))^{1/2}}f_{4h}(F(u,v);\mu)\hspace{3cm}\]
\[\hspace{3cm}{}\ll_F \frac{B}{N(\mu)}
\exp\left\{\sum_{2<p\le B}\frac{\rho(p)(f_{4h}(p;\mu)-1)}{p}\right\}.\]
However
\[\sum_{2<p\le B}\frac{\rho(p)(f_{4h}(p;\mu)-1)}{p}=
\sum_{2<p\le B}\frac{\rho(p)(f_{4h}(p)-1)}{p}
+O\left(\sum_{p\mid N(\mu)}\frac{1}{p}\right).\]
The error term is $O(\log\log\log N(\mu))$, by the argument used for
(\ref{lll}), whence
\[\sum_{\max(|u|,|v|)\le (B/N(\mu))^{1/2}}f_{4h}(F(u,v);\mu)\hspace{3cm}\]
\[\hspace{2cm}{}\ll_F \frac{B(\log\log N(\mu))^A}{N(\mu)}
\exp\left\{\sum_{2<p\le B}\frac{\rho(p)(f_{4h}(p)-1)}{p}\right\}\]
for a suitable constant $A$, so that (\ref{yi}) yields
\[\sum_{t^{1/2}\max(|u|,|v|)\le B^{1/2}}
  \left|\sum_{\substack{x^2+y^2=t^2F(u,v)\\ \gcd(t,x,y)=1}}
  \left(\frac{x+iy}{|x+iy|}\right)^{4h}\right|\hspace{3cm}\]
  \[\hspace{1cm}{}\ll_F B(\log\log B)^A(\log B)
 \exp\left\{\sum_{2<p\le B}\frac{\rho(p)(f_{4h}(p)-1)}{p}\right\}.\]
We now apply the inequality (\ref{34}) to estimate $f_{4h}(p)$ when
$p\equiv 1($mod 4), and deduce that
\[\sum_{t^{1/2}\max(|u|,|v|)\le B^{1/2}}
  \left|\sum_{\substack{x^2+y^2=t^2F(u,v)\\ \gcd(t,x,y)=1}}
  \left(\frac{x+iy}{|x+iy|}\right)^{4h}\right| \hspace{3cm}\]
  \beql{123}
  \hspace{2cm}{}\ll_F B(\log\log B)^A(\log B)
  \exp\left\{-S_1+\tfrac32 S_2+\tfrac14 S_3\right\},
  \eeq
 with
 \[S_1=\sum_{2<p\le B}\frac{\rho(p)}{p}\]
 \[S_2=\sum_{\substack{2<p\le B\\ p\equiv 1(\mathrm{mod}\, 4)}}\frac{\rho(p)}{p}\]
 and
 \[S_3= \sum_{2<p\le B}\frac{\rho(p)g_{8h}(p)}{p}.\]
 
We observed above that $\rho(p)$ is just the number of solutions of
$F(X,1)$ modulo $p$, if $p$ is large enough in terms of $F$. Indeed,
if we write $K=\QQ(\theta)$ where $\theta$ is a root of $F(X,1)$, we
see from Dedekind's Theorem that $\rho(p)$ is the number of first
degree prime ideals $P$ of $K$ above $p$, at least if $p$ is large enough
in terms of $F$.  It follows that
\[S_1=\sum_{2<N(P)\le B}\frac{1}{N(P)}+O_F(1),\]
since prime ideals of degree 2 or more contribute $O_F(1)$. The Prime
Ideal Theorem then shows that
\beql{E1}
S_1=\log\log B+O_F(1).
\eeq
In order to examine $S_2$ and $S_3$ we will consider the factorization
of rational primes in the field $L=K(i)$. This must be a quadratic
extension of $K$, since if $\QQ(i)\subset K$ it would follow that
$\theta$ is quadratic over $\QQ(i)$, which is impossible since $F(X,1)$
was assumed to be irreducible over $\QQ(i)$. When $p\equiv 1($mod 4)
is a norm $N(\pi)$ over $\QQ(i)$ one sees that each prime ideal $P$ of
$K$ lying above $p$ will split as $(P,\pi)(P,\overline{\pi})$ over
$L$, so that $\rho(p)$ is half the number of first degree prime ideals
of $L$ lying over $p$.  On the other hand, there cannot be a first degree
prime ideal $P$ of $L$ lying above a rational prime $p\equiv 3($mod 4), 
since then $N_{L/\QQ(i)}(P)$ would be a first degree
prime ideal of $\QQ(i)$ above $p$, which is impossible. It follows
that
\[S_2=\frac12\sum_{2<p\le B}\frac{\rho_L(p)}{p},\]
where $\rho_L(p)$ counts first degree primes of $L$ above $p$. The
Prime Ideal theorem for $L$ then yields
\beql{E2}
S_2=\tfrac12\log\log B+O_F(1).
\eeq
To handle $S_3$ we define a non-trivial Hecke Gr\"{o}ssencharacter on
prime ideals of $L$ by setting $\chi_h(P)=(\pi/|\pi|)^{4h}$ if
$N_{L/\QQ(i)}(P)=(\pi)$. With this definition we will have
\[\rho(p)g_{8h}(p)=\sum_{N_{L/\QQ}(P)=p}\chi_{2h}(p).\]
In analogy with Lemma \ref{pntg} we deduce from
the Prime Ideal Theorem with Gr\"{o}ssencharacter that
\[\sum_{N(A)\le x}\Lambda(A)\chi_h(A)\ll_F 
  x\exp\left\{-c\frac{\log x}{\sqrt{\log x}+\log h}\right\}(\log xh)^4\]
  uniformly for $x\ge 2$ and $k\in\NN$, where $A$ runs over non-zero
  integral ideals of $L$. Following a similar argument to that used in
  \S \ref{pt2}  we now deduce that
\[S_3=\sum_{3\le p\le N}\frac{\rho(p)g_{8h}(p)}{p}\ll_F \log\log h.\]
We can now feed this estimate, along with (\ref{E1}) and (\ref{E2}),
into (\ref{123}) to deduce that
\[\sum_{t^{1/2}\max(|u|,|v|)\le B^{1/2}}
  \left|\sum_{\substack{x^2+y^2=t^2F(u,v)\\ \gcd(t,x,y)=1}}
  \left(\frac{x+iy}{|x+iy|}\right)^{4h}\right|\hspace{3cm}\]
\[\hspace{3cm}{}  \ll_F B(\log\log B)^A(\log h)^A(\log B)^{3/4},\]
for a suitable numerical constant $A$, and Lemma \ref{lt2} follows.

\end{document}